\documentclass[11pt, twoside]{article} 
\usepackage{amsmath}
\usepackage{amssymb}
\usepackage{amsthm}
\usepackage{amsfonts}
\usepackage{graphicx} 
\usepackage[T1]{fontenc}
\usepackage[utf8]{inputenc}  
\usepackage{lineno,hyperref}
\usepackage[top=4cm,bottom=4.5cm,left=3cm,right=3cm,asymmetric]{geometry}
\usepackage{epstopdf}
\usepackage{indentfirst}
\usepackage{anyfontsize}
\usepackage{atbegshi}
\usepackage{tikz}
\usepackage{float}

\newtheorem{theorem}{Theorem}

\usepackage{natbib}
\setcitestyle{numbers,sort&compress}
\providecommand{\keywords}{\vspace{0.8em}\textbf{Keywords: }}

\title{ \large \textbf{\uppercase{A note on spanning trees with a specified degree sequence}}}

\author{Mar\'ia Elena Mart\'inez-Cuero \\ {\small\it Departamento de Matem\'aticas }\\\ \textit{\small Universidad Aut\'onoma Metropolitana-Iztapalapa}  \\ {\small\tt sherlyroses@hotmail.com}\\ \\ Eduardo Rivera-Campo \\ {\small\it Departamento de Matem\'aticas }\\\ \textit{\small Universidad Aut\'onoma Metropolitana-Iztapalapa} \\  {\small\tt erc@xanum.uam.mx} }

\date{}

\renewenvironment{abstract}
               {\small\list{}{\rightmargin\leftmargin}%
                \item[\hspace{10mm}\normalfont\textsc{abstract}.]\relax}
               {\endlist}

\usepackage{fancyhdr}
\usepackage{sectsty}
 
\allsectionsfont{\fontsize{12}{12}\selectfont\centering\mdseries\scshape}

\newcommand\blfootnote[1]{%
  \begingroup
  \renewcommand\thefootnote{}\footnote{#1}%
  \addtocounter{footnote}{-1}%
  \endgroup
}

\setcitestyle{square}
\makeatletter
\long\def\@makecaption#1#2{
\vskip\abovecaptionskip
\sbox\@tempboxa{#1. #2}
\ifdim \wd\@tempboxa >\hsize
#1. #2\par
\else
\global \@minipagefalse
\hb@xt@\hsize{\hfil\box\@tempboxa\hfil}
\fi
\vskip\belowcaptionskip}
\makeatother

\newcommand{\QED}{\hfill {\qed}}

\begin{document}
\maketitle

\begin{abstract} 
We give an Ore-Type condition sufficient for a graph $G$ to have a spanning tree with a specified degree sequence.

\end{abstract}
\blfootnote{Partially supported by Conacyt, México.}

\section{Introduction} 

\par  O. Ore \cite{O60} proved that if $G$ is a graph with $n$ vertices such that $d(u) + d(v) \geq n - 1$ for each pair $u, v$ of non-adjacent vertices, then $G$ contains a hamiltonian path. This result has been generalized in many directions.

S. Win \cite{W75} showed that if $r \geq 2$ is an integer and $G$ is a connected graph with $n$ vertices such that $d(u_1) + d(u_2) + \cdots + d(u_r) \geq n-1$ for each set of $r$ independent vertices of $G$, then $G$ has a spanning tree with maximum degree at most $r$.

Years later, H. Broersma and H. Tuinstra \cite{BT98} showed that if $s \geq 2$ is an integer and $G$ is a connected graph with $n$ vertices such that $d(u) + d(v) \geq n-s+1$ for
each pair $u, v$ of non-adjacent vertices, then $G$ contains a spanning tree with at most $s$ vertices with degree 1.

E. Rivera-Campo \cite{RC12} gave a condition on the graph $G$ that bounds the degree of each vertex in a certain spanning tree $T$ of $G$  and the number of vertices of $T$ with degree 1.

\begin{theorem}

Let $n$, $k$, and $d_1, d_2, \ldots, d_n$ be integers with $1 \leq k \leq n-1$ and $2 \leq d_1 \leq d_2 \leq  \cdots \leq  d_n \leq n-1$. If $G$ is a $k$-connected graph with vertex set $V(G) = \left\{ w_1, w_2, \ldots, w_n \right\}$ such that $d(u) + d(v) \geq n - 1 - \sum_{j=1}^{n} {(d_j - 2)}$ for each pair $u, v$ of non-adjacent vertices, then $G$ contains a spanning tree $T$ with at most $2 + \sum_{j=1}^{n} {(d_j - 2)} $ vertices with degree 1 and such that $d_T(w_j) \leq d_j$ for $j=1, 2,\ldots, n$.

\end{theorem}

Let $n$ be a positive integer. An \emph {arboreal sequence} is a sequence of positive integers $d_1, d_2, \ldots, d_n$ such that $\sum_{j=1}^{n} {d_j} =2(n-1) $. It is well known that a sequence $\sigma=d_1, d_2, \ldots, d_n$ is arboreal if and only if there is a tree whose vertices have degrees $d_1, d_2, \ldots, d_n$.

Let $\sigma = d_1, d_2, \ldots, d_n$ be an arboreal sequence and $G$ be a labelled graph with vertex set $V(G) = \left\{w_1, w_2, \ldots, w_n\right\}$. A spanning tree $T$ of $G$ has degree sequence $\sigma$ if $d_T(w_i)=d_i$ for $i=1,2,\ldots, n$. In this note we prove the following result:

\begin{theorem}
\label{principal}
Let  $n \geq 4$ be an integer and $G$ be a labelled graph with vertex set $V(G) = \left\{w_1, w_2, \ldots, w_n\right\}$. If $d(u) + d(v) \geq \frac {3n-1} {2}$ for each pair $u, v$ of non-adjacent vertices, then $G$ contains a spanning tree $T$ with degree sequence $\sigma$ for each arboreal sequence $\sigma = d_1, d_2, \ldots, d_n$  with $1\leq d_i \leq 3$ for $i=1,2,\ldots, n$.
\end{theorem}

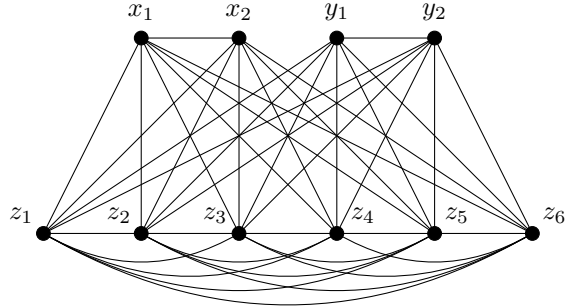
\begin{figure}[!ht]
 
\begin{center}
\begin{tikzpicture}[scale=1.3]

 \filldraw[color=black]  (0,0) circle (2pt) node[above left] { \small {$z_1$}}; 
 \filldraw[color=black]  (1,0) circle (2pt) node[above left] {\small {$z_2$}}; 
 \filldraw[color=black] (2,0) circle (2pt) node[above left, xshift=-.2mm] {\small {$z_3$}}; 
 \filldraw[color=black] (3,0) circle (2pt) node[above right, xshift=.5mm] {\small {$z_4$}}; 
 \filldraw[color=black] (4,0) circle (2pt) node[above right] {\small {$z_5$}}; 
 \filldraw[color=black] (5,0) circle (2pt) node[above right] {\small{$z_6$}}; 
 \filldraw[color=black] (1,2) circle (2pt) node[above,yshift=1mm ] { \small {$x_1$}}; 
 \filldraw[color=black] (2,2) circle (2pt) node[above,yshift=1mm ] { \small {$x_2$}}; 
 \filldraw[color=black] (3,2) circle (2pt) node[above,yshift=1mm ]{ \small {$y_1$}}; 
 \filldraw[color=black] (4,2) circle (2pt) node[above,yshift=1mm ]{ \small {$y_2$}}; 
 
\draw[thin,-] (0,0) -- (5,0); 
\draw[thin,-] (0,0) to[bend right] (2,0);
\draw[thin,-] (0,0)to[bend right]  (3,0);
\draw[thin,-] (0,0) to[bend right](4,0);
\draw[thin,-] (0,0) to[bend right] (5,0);
\draw[thin,-] (1,0)to[bend right]  (3,0); 
\draw[thin,-] (1,0)to[bend right]  (4,0);  
\draw[thin,-] (1,0)to[bend right]  (5,0); 


\draw[thin,-] (2,0)to[bend right]  (4,0);  
\draw[thin,-] (2,0)to[bend right]  (5,0);

\draw[thin,-] (3,0)to[bend right]  (5,0);

\draw[thin,-] (1,2) -- (2,2); 
\draw[thin,-] (1,2) -- (0,0); 
\draw[thin,-] (1,2) -- (1,0);
\draw[thin,-] (1,2) -- (2,0);
\draw[thin,-] (1,2) -- (3,0);
\draw[thin,-] (1,2) -- (4,0);
\draw[thin,-] (1,2) -- (5,0);

\draw[thin,-] (2,2) -- (0,0); 
\draw[thin,-] (2,2) -- (1,0);
\draw[thin,-] (2,2) -- (2,0);
\draw[thin,-] (2,2) -- (3,0);
\draw[thin,-] (2,2) -- (4,0);
\draw[thin,-] (2,2) -- (5,0);

\draw[thin,-] (3,2) -- (4,2); 
\draw[thin,-] (3,2) -- (0,0); 
\draw[thin,-] (3,2) -- (1,0);
\draw[thin,-] (3,2) -- (2,0);
\draw[thin,-] (3,2) -- (3,0);
\draw[thin,-] (3,2) -- (4,0);
\draw[thin,-] (3,2) -- (5,0);

\draw[thin,-] (4,2) -- (0,0); 
\draw[thin,-] (4,2) -- (1,0);
\draw[thin,-] (4,2) -- (2,0);
\draw[thin,-] (4,2) -- (3,0);
\draw[thin,-] (4,2) -- (4,0);
\draw[thin,-] (4,2) -- (5,0);

\end{tikzpicture}

\end{center} 

\caption{Graph $G_2$}  \label{Fig 1}
\end{figure}

For each positive integer $k$ let $X_k = \{ x_1, x_2, \ldots, x_k\}$, $Y_k = \{ y_1, y_2, \ldots, y_k\}$ and $Z_k = \{ z_1, z_2, \ldots, z_{2k+2}\}$ be pairwise disjoint sets of vertices and let $G_k$ be the complete graph with vertex set $X_k \cup Y_k \cup Z_k$ with all edges $x_iy_i$, $1\leq i, j\leq k$ removed. See Fig. \ref{Fig 1}  for the case $k=2$.

We claim that the graph $G_k$ contains no spanning tree $T$ such that $d_T(x_i) = d_T(y_i) = 3$ for $i = 1, 2, \ldots, k$ and $d_T(z_j) = 1$ for $j=1, 2, \ldots, 2k+2$; for if $T$ is such a tree, then $T - Z_k$ would be a spanning tree of the subgraph $G \left[ X_k \cup Y_k \right]$ of $G$, induced by the set $X_k \cup Y_k$,  which is not possible since $G \left[ X_k \cup Y_k \right]$ is not connected. On the other hand, if $u$ and $v$ are non-adjacent vertices of $G_k$, without loss of generality we may assume $u \in X_k$ and $v \in Y_k$. Therefore 

$$d_{G_k}(u) + d_{G_k}(v) = 2((k-1) +(2k+2)) = 6k+2 = \frac {3n-2} {2}$$ 

 where $n= 4k+2$ is the number of vertices of $G_k$. This shows that the degree-sum condition in Theorem \ref{principal} is tight.

Whenever possible we follow the notation of J. A. Bondy and U. S. R. Murty \cite{BM}.

\section{Proof of Theorem \ref{principal}}

Suppose the result is false. Then for certain integer $n \ge 4$ and certain arboreal sequence $\sigma = d_1, d_2, \ldots, d_n$  with $1\leq d_i \leq 3$ for $i=1,2,\ldots, n$ there exists a counterexample. That is a labelled graph $G$ with vertex set $V(G) = \left\{w_1, w_2, \ldots, w_n\right\}$ such that $G$ contains no spanning tree with degree sequence $\sigma$ while  $d(u) + d(v) \geq \frac {3n-1} {2}$ for each pair $u, v$ of non-adjacent vertices of $G$. We choose $G$ with the maximum possible number of edges while remaining a counter example with $n$ vertices.

Since $\sigma$ is an arboreal sequence of order $n$, a counterexample cannot be a complete graph of order $n$. Let $u, v$ be non-adjacent vertices of $G$. By the choice of $G$, the graph $G+uv$ is not a counterexample and contains a spanning tree $T$ with degree sequence $\sigma$. Therefore $G$ contains a spanning forest $F=T-uv$ with exactly two components $T_u$ and $T_v$ with $u\in V(T_u)$ and $v\in V(T_v)$ such that $d_F(u) = d_i -1$, $d_F(v)=d_j -1$ and $d_f(w_r) =d_r$ for each vertex $u \neq w_r \neq v$, where $i$ and $j$ are such that $u = w_i$ and $v = w_j$.

Orient the edges of $F$ in such a way that $T_u$ and $T_v$ become outdirected trees  $\overrightarrow{T_u}$ and $\overrightarrow{T_v}$ with roots $u$ and $v$, respectively. For each vertex $u \neq w \neq v$ let $w^-$ denote the unique vertex of $G$ such that the edge $w^- w$ is oriented from $w^-$ to $w$ in $\overrightarrow{F}$. Let

$$A_u =\{ y^{-} \in  V(T_u): uy \in E(G)\}, B_u =\{x \in  V(T_u): vx \in E(G)\},$$  $$A_v =\{ y^{-} \in  V(T_v): vy \in E(G)\}
 \text{ and } B_v =\{x \in  V(T_v): ux \in E(G)\}.$$ \negthickspace

\begin{figure}[!ht]
\begin{center}
\begin{tikzpicture}[scale=0.9]

  \filldraw[color=black] (3,0) circle (2pt); 
  \filldraw[color=black] (1,1) circle (2pt); 
  \filldraw[color=black] (2,1) circle (2pt); 
  \filldraw[color=black] (3,1) circle (2pt); 
  \filldraw[color=black] (6,1) circle (2pt);
  \filldraw[color=black] (7,1) circle (2pt);
  \filldraw[color=black] (4,2) circle (2pt) node[above] {\small {$u$}}; 
  \filldraw[color=black] (5,2) circle (2pt) node[above] {\small {$v$}}; 
  \filldraw[color=black] (2,3) circle (2pt);
  \filldraw[color=black] (3,3) circle (2pt);
  \filldraw[color=black] (6,3) circle (2pt);
  \filldraw[color=black] (7,3) circle (2pt);
  \filldraw[color=black] (3,4) circle (2pt);
  \filldraw[color=black] (6,4) circle (2pt);
  \filldraw[color=black] (1,2.5) circle (2pt);
  \filldraw[color=black] (1,3.5) circle (2pt);

\draw[->, >=latex] (3,1) to (2.05,1);
\draw[->, >=latex] (2,1) to (1.05,1);
\draw[->, >=latex] (3,1) to (3,0.05);
\draw[->, >=latex] (4,2) to (3,1.05);
\draw[->, >=latex] (4,2) to (3.08,2.98);
\draw[->, >=latex] (3,3) to (2.05,3);
\draw[->, >=latex] (2,3) to (1.07,3.5);
\draw[->, >=latex] (2,3) to (1.07,2.5);
\draw[->, >=latex] (3,3) to (3,3.93);

\draw[->, >=latex] (5,2) to (5.97,2.97);
\draw[->, >=latex] (5,2) to (5.99,1.07);
\draw[->, >=latex] (6,3) to (6.95,3);
\draw[->, >=latex] (6,3) to (6,3.95);
\draw[->, >=latex] (6,1) to (6.95,1);


 \filldraw[color=black] (11,0) circle (2pt); 
 \filldraw[color=black] (9,1) circle (2pt); 
 \filldraw[color=black] (10,1) circle (2pt); 
 \filldraw[color=black] (11,1) circle (2pt);
 \filldraw[color=black] (14,1) circle (2pt); 
 \filldraw[color=black] (15,1) circle (2pt);
 \filldraw[color=black] (12,2) circle (2pt) node[above] {\small {$u$}};
 \filldraw[color=black] (13,2) circle (2pt) node[above] {\small {$v$}};
 \filldraw[color=black] (10,3) circle (2pt);
 \filldraw[color=black] (11,3) circle (2pt);
 \filldraw[color=black] (14,3) circle (2pt);
 \filldraw[color=black] (15,3) circle (2pt);
 \filldraw[color=black] (11,4) circle (2pt);
 \filldraw[color=black] (11,4) circle (2pt);
 \filldraw[color=black] (14,4) circle (2pt);
 \filldraw[color=black] (9,2.5) circle (2pt);
 \filldraw[color=black] (9,3.5) circle (2pt);
 
\draw [-] (11,0) to (11,1);
\draw [-] (11,1) to (9,1);
\draw [-] (11,1) to (12,2);
\draw[-] (12,2) to (11,3); 
\draw[-] (12,2)to [bend left]  (10,3); 
\draw [-] (10,3) to (9,2.5);
\draw [-] (10,3) to (9,3.5);
\draw [-] (11,3) to (11,4);
\draw[-] (11,3) to [bend left] (13,2);
\draw[-] (13,2) to (14,1);
\draw[-] (13,2) to (14,3);
\draw[-] (14,3) to (14,4);
\draw[-] (14,3) to (15,3);
\draw[-] (14,1) to (15,1);

\end{tikzpicture}
\end{center}
\caption{Forest $F$ with $z^- \in A_u \cap B_u$ (left). Tree $(F-z^- z) + \{ uz,vz^-  \}$(right).}  \label{Fig2}
\end{figure}
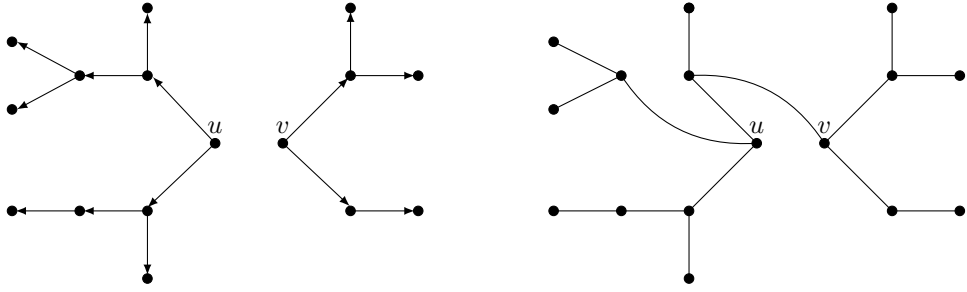

Notice that $|A_u \cap B_u|=0$, for if $z^{-}\in A_u \cap B_u$, then $(F-z^{-}z)+ \{uz, vz^- \}$ would be a spanning tree of $G$ with degree sequence $\sigma$, which is not possible (See Fig. \ref{Fig2}). Analogously $|A_v \cap B_v|=0$ and therefore $$|A_u|+|B_u|=|A_u \cup B_u|\leq n_u \text{ and }|A_v|+|B_v|=|A_v \cup B_v|\leq n_v,$$ where $n_u$ and $n_v$ are the number of vertices of  $T_u$ and $T_v$, respectively. 

In an abuse of notation, for each vertex $x$ of $G$ we denote by $d_u(x)$ and $d_v(x)$ the number of vertices of $T_u$ and $T_v$, respectively, which are adjacent to $x$ in $G$. Clearly $$|B_u|=d_{u}(v) \text{ and } |B_v|= d_{v}(u).$$ Also notice that $$|A_u|\geq \frac{d_u (u)}{2} \text{ and } |A_v|\geq \frac{d_v (v)}{2}$$

since the out-degree of each vertex of $\overrightarrow{F}$ is at most 2. Then $$\frac{d_u (u)}{2}+ d_u (v) \leq |A_u|+|B_u| \leq n_u  \text{ and }
\frac {d_v (v)}{2}+ d_v (u) \leq |A_v|+|B_v| \leq n_v .$$ Therefore  $$d_u(u)+2d_u(v) \leq 2 n_u \text{ and } d_v(v)+2d_v (u) \leq 2n_v.$$ Since $u$ is not  adjacent to $u$ and $v$ is not adjacent to $v$ in $G$, $d_u(u) \leq n_u -1$ and $d_v(v) \leq n_v -1$. Adding these to the previous inequalities we obtain $$2d_u(u)+2d_v(v) \leq 3 n_u -1 \text{ and } 2d_v(v)+2d_v (u) \leq 3n_v -1.$$ These imply $$2d(u)+2d(v)=(2d_u (u)+2d_v(u))+(2d_v(v)+2d_u(v))\leq 3(n_u + n_v)-2 = 3n-2,$$ which is not possible since $d(u)+ d(v)\geq \frac{3n-1}{2}.$ \QED

\section{FINAL REMARKS}

With the same approach, we can prove the following generalization of Theorem \ref{principal}.

\begin{theorem}
\label{general}
Let $r\geq 2$ be an integer and $G$ be a labelled graph with vertex set $V(G) = \left\{ v_1, v_2, \ldots, v_n \right\}$ with $n\geq r+1$. If $d(u)+d(v)\geq\frac{(2r-3)n-(2r-5)}{r-1}$ for each pair $u$, $v$ of non-adjacent vertices, then $G$ has a spanning tree $T$ with degree sequence $\sigma$ for each arboreal sequence $\sigma=d_1, d_2, \ldots, d_n $ with $1\leq d_i\leq r$ for $i=1, 2, \ldots, n$.

\end{theorem}

Let $r\geq 2$ be an integer. For each positive integer $k$ let $X_k=\left\{x_1, x_2, \ldots, x_k \right\}, Y_k= \left\{ y_1, y_2, \ldots, y_k \right\}$ and $Z_{k,r}= \left\{ z_1, z_2, \ldots, z_{2k(r-2)+2} \right\}$ be pairwise disjoint vertex sets and let $G_{k,r}$ be the complete graph with vertex set 
$X_k \cup Y_k\cup Z_{k,r}$ with all edges $x_i y_i, 1\leq i,j\leq k$ removed. As for the graphs $G_k=G_{k,3}$ in the introduction, we claim that the graphs $G_{k,r}$ show that the degree-sum condition in Theorem \ref{general}  is also tight.

\maketitle

\end{document}